\input amstex
\documentstyle{amsppt}
\document

\magnification 1100

\def\eps{{\varepsilon}}

\def\1b{{\bold 1}}

\def\hb{{\bold h}}

\def\ul{{\bold u}}

\def\xb{{\bold x}}

\def\Ab{{\bold A}}

\def\Ub{{\bold U}}

\def\Rep{\text{Rep}\,}

\def\res{\text{res}\,}

\def\End{\text{End}\,}
\def\Lin{\text{Lin}\,}

\def\Ker{\text{Ker}\,}

\def\Im{\text{Im}\,}

\def\CC{{\Bbb C}}

\def\LL{{\Bbb L}}

\def\NN{{\Bbb N}}
\def\PP{{\Bbb P}}

\def\ZZ{{\Bbb Z}}

\def\Cc{{\Cal C}}

\def\Ec{{\Cal E}}
\def\Fc{{\Cal F}}

\def\Hc{{\Cal H}}

\def\Lc{{\Cal L}}

\def\Oc{{\Cal O}}

\def\Uc{{\Cal U}}
\def\Vc{{\Cal V}}
\def\Wc{{\Cal W}}

\def\and{{\quad\text{and}\quad}}
\def\where{{\quad\text{where}\quad}}
\def\st{{\, \mid\,}}
\def\ds{\displaystyle}

\def\sss{\scriptscriptstyle}                
\def\qed{\hfill $\sqcap \hskip-6.5pt \sqcup$}        
\overfullrule=0pt                                    

\def\ua{{\alpha}}
\def\ub{{\beta}}
\def\ul{{\lambda}}
\def\um{{\mu}}
\def\un{{\nu}}
\def\uo{{\omega}}

\def\var{\varepsilon}
\def\heck{{\buildrel{\sss k}\over\to}}
\def\hecl{{\buildrel{\sss l}\over\to}}

\newdimen\Squaresize\Squaresize=14pt
\newdimen\Thickness\Thickness=0.5pt
\def\Square#1{\hbox{\vrule width\Thickness
	      \vbox to \Squaresize{\hrule height \Thickness\vss
	      \hbox to \Squaresize{\hss#1\hss}
	      \vss\hrule height\Thickness}
	      \unskip\vrule width \Thickness}
	      \kern-\Thickness}
\def\Vsquare#1{\vbox{\Square{$#1$}}\kern-\Thickness}

\centerline{\bf ON THE K-THEORY OF THE}
\centerline{\bf CYCLIC QUIVER VARIETY} 

\vskip 15mm

\centerline{\bf M. VARAGNOLO and E. VASSEROT
\footnote{
Both authors are partially supported by EEC grant
no. ERB FMRX-CT97-0100.\hfill\break
}}

\vskip3cm

\centerline{\bf 1. Introduction and notations.}

\vskip3mm

In \cite{N1} Nakajima has defined a new class of varieties, called 
{\bf quiver varieties}, associated to any quiver.
In the particular case of the Dynkin quiver of finite type $A$, 
these varieties are related to partial flag manifolds.
For any quiver with no edge loops, he proved in \cite{N2} that the 
corresponding
simply laced Kac-Moody algebra acts on the top homology groups
via a convolution product.
The purpose of this paper is to compute the convolution product 
on the equivariant K-groups of the cyclic quiver variety, 
generalizing the previous works \cite{GV},\cite{V}. 
It is expected that equivariant $K$-groups should
give an affinization of the quantized enveloping algebra
of the corresponding Kac-Moody algebra.
Similar algebras, called {\bf toroidal algebras}, have already been studied 
in the cyclic case (e.g. \cite{GKV}, \cite{VV}).
Surprisingly the convolution operators we get do not satisfy exactly the 
relations of the quantized toroidal algebra.
We obtain a twisted version, denoted by $\Ub_{q,t}$, of the 
latter. Fix an integer $n\geq 3$ and two parameters $q,t\in\CC^\times$. 
Then, $\Ub_{q,t}$ is the complex unital associative 
algebra generated by $\xb^{^+}_{i,k},$ $\xb^{^-}_{i,k},$ 
$\hb^{^\pm}_{i,\pm l},$ where $i\in\ZZ/n\ZZ,$ $k\in\ZZ$, $l\in\NN$.
The relations are expressed in term of the formal series
$$\xb^{^+}_i(z)=\sum_{k\in\ZZ}{\xb^{^+}_{i,k}\cdot z^{^{-k}}},\quad
\xb^{^-}_i(z)=\sum_{k\in\ZZ}{\xb^{^-}_{i,k}\cdot z^{^{-k}}}\and
\hb_i^{^\pm}(z)=\sum_{k\geq 0}{\hb^{^\pm}_{i,\pm k}\cdot z^{^{\mp k}}},$$
as follows

$$\hb_{i,0}^{\eps}\,\hb_{i,0}^{-\eps}=1,\qquad
[\hb^{\eps}_i(z),\hb^{\tau}_j(w)]=0$$
$$(1-q^{-1}t^{-1})\,[\xb^{+}_i(z),\xb^{-}_j(w)]=\delta_{ij}\,
\epsilon(z/w)\bigl(\hb^{+}_i(w)-\hb^{-}_i(z)\bigr)$$
$$(z-qtw)^{\tau}\hb^{\eps}_i(z)\,\xb^{\tau}_i(w)=
(qtz-w)^{\tau}\xb^{\tau}_i(w)\,\hb^{\eps}_i(z)$$
$$(tz-w)^{\tau}\hb^{\eps}_{i+1}(z)\,\xb^{\tau}_i(w)=
(z-qw)^{\tau}\xb^{\tau}_i(w)\,\hb^{\eps}_{i+1}(z)$$
$$(z^{\eps}-qtw^{\eps})\,\xb^{\eps}_i(z)\,\xb^{\eps}_i(w)=
(qtz^{\eps}-w^{\eps})\,\xb^{\eps}_i(w)\,\xb^{\eps}_i(z)$$
$$(tz^{\eps}-w^{\eps})\,\xb^{\eps}_{i+1}(z)\,\xb^{\eps}_i(w)=
(z^{\eps}-qw^{\eps})\,\xb^{\eps}_i(w)\,\xb^{\eps}_{i+1}(z)$$
$$t^{\tau}\,\xb^{\eps}_i(z_{_1})\,\xb^{\eps}_i(z_{_2})\,\xb^{\eps}_{i+\tau}(w)
-\bigl(q^{\tau}t^{\tau}+1\bigr)\xb^{\eps}_i(z_{_1})
\,\xb^{\eps}_{i+\tau}(w)\,\xb^{\eps}_i(z_{_{2}})+$$
$$+q^{\tau}\,\xb^{\eps}_{i+\tau}(w)\,\xb^{\eps}_i(z_{_{1}})\,
\xb^{\eps}_i(z_{_2})+\{z_{_1}\leftrightarrow z_{_2}\}=0$$
$$[\xb^{\eps}_i(z),\xb^{\eps}_j(w)]=0\qquad\text{if}\quad i\neq j,j\pm 1,$$

\vskip2mm

\noindent where $\epsilon(z)=\sum_{n\in\ZZ}z^n$ and $\eps,\tau\in\{+1,-1\}$.
Observe that we recover the usual Drinfeld relations after the specialization 
$q=t$. Our work relies very much on \cite{N2},\cite{N3}. 
In Section 2 we make several recollections. 
Theorem 1 is not stated explicitely in \cite{N2} but Nakajima 
told us it was known to him. We mention it here since it extends and 
clarifies some result in \cite{L}. The convolution product is computed in 
Sections 3 and 4. Another reason to give both Theorem 1 and Theorem 2
is the general philosophy mentioned in \cite{GKV} and in Nakajima's works
that convolution algebras associated to moduli spaces of torsion free
sheaves on surfaces should lead to affinization of affine algebras
and their quantum analogues. 

\vskip1cm

\centerline{\bf 2. Quiver varieties.} 

\vskip3mm

\noindent{\bf 2.1.} Let us recall a few definitions from \cite{N2}.
Fix a finite subgroup $\Gamma\subset SL_2(\CC).$
Let $\Rep(\Gamma)$ be the set of isomorphism classes of finite dimensional
representations of $\Gamma$ and let $L\in\Rep(\Gamma)$ be the obvious
2-dimensional representation. 
Let $X(\Gamma)=\{S_k\st k=1,2,...,n\}$ be the set of the simple 
representations of $\Gamma$. 
Given finite dimensional vector spaces $V,W$ set
$$M_{W,V}=\Lin(V,L\otimes V)\oplus\Lin(W,\wedge^2L\otimes V)\oplus\Lin(V,W).$$
The group $GL(V)$ acts on $M_{W,V}$ as follows : for any $g\in GL(V)$,
$$g\cdot(B,i,j)=\bigl((1\otimes g)\circ B\circ g^{-1},(1\otimes g)\circ i, 
j\circ g^{-1}\bigr).$$
Let us consider the map
$$\mu_{W,V}\,:\quad M_{W,V}\to\Lin(V,\wedge^2L\otimes V),\quad
(B,i,j)\mapsto B\wedge B-i\circ j,$$
where $B\wedge B$ is the projection of 
$B\circ B\in\Lin(V,L^{\otimes 2}\otimes V)$ to $\Lin(V,\wedge^2L\otimes V)$.
A triple $(B,i,j)\in\mu_{W,V}^{-1}(0)$ is {\bf stable}
if and only if there is no nontrivial
subspace $V'\subseteq\Ker j$ such that $B(V')\subseteq L\otimes V'$.
Let $\mu_{W,V}^{-1}(0)^s$ be the set of stable triples.
The group $GL(V)$ acts freely on $\mu_{W,V}^{-1}(0)^s$ 
(see \cite{N2, Lemma 3.10}).
Suppose now that $W,V$ are $\Gamma$-modules.
Then $M_{W,V}$ is endowed with an obvious $\Gamma$-action. 
Let $\mu_{W,V}^{-1}(0)^{s,\Gamma}\subseteq\mu_{W,V}^{-1}(0)^s$
be the fixpoint variety. 
Let $GL_\Gamma(V)$ and $V_k=\Lin_\Gamma(S_k,V)$ be the $\Gamma$-invariants.
Let 
$$\Lambda_{W,V}=\mu_{W,V}^{-1}(0)^{s,\Gamma}/GL_\Gamma(V)\and
N_{W,V}=\mu^{-1}_{W,V}(0)^\Gamma//GL_\Gamma(V)$$
be the Nakajima varieties.
By construction, $\Lambda_{W,V}$ is a smooth variety (or the empty set).
Put $\Lambda_W=\bigcup_V\Lambda_{W,V}$ and $N_W=\bigcup_VN_{W,V}$.
Consider the morphism $\pi\,:\,\Lambda_{W,V}\to N_{W,V}$ induced by the 
chain of maps 
$$\mu^{-1}_{W,V}(0)^{s,\Gamma}\hookrightarrow
\mu^{-1}_{W,V}(0)^\Gamma\to\mu^{-1}_{W,V}(0)^\Gamma//GL_\Gamma(V).$$

\vskip3mm

\noindent{\bf 2.2.} 
Let us now recall a few results from \cite{N2},\cite{L}.
A triple $(B,i,j)\in\mu^{-1}_{W,V}(0)$ is {\bf costable} if and only if 
there is no subspace $V'\subseteq V$ such that
$$B(V')\subseteq L\otimes V',\qquad\Im i\subseteq\wedge^2L\otimes V'
\and V'\neq V.$$
Let $\mu^{-1}_{W,V}(0)^c$ be the set of costable triples and put
$$\mu^{-1}_{W,V}(0)^{s,c}=\mu^{-1}_{W,V}(0)^c\cap\mu^{-1}_{W,V}(0)^s.$$
For any triple $\xi=(B,i,j)\in\Lambda_{W,V}$ the {\bf dual},
$\xi^*$, is the triple
$(B^*,j^*,i^*)\in\Lambda_{W^*,V^*}$, where
$$B^*=\wedge^2L\otimes B^t,\quad j^*=\wedge^2L\otimes j^t,\quad
i^*=\wedge^2L\otimes i^t,$$
and ${}^t$ stands for the transposed map.
Since $W^*\simeq W$ and $V^*\simeq V$ as $\Gamma$-modules, we will view
$\xi^*$ as an element of $\Lambda_{W,V}.$
Then, $\xi\in\mu_{W,V}^{-1}(0)^s$ if and only if $\xi^*\in\mu_{W,V}^{-1}(0)^c$.

\vskip3mm

\noindent{\smc Lemma 1.} {\it Fix $\xi\in\mu^{-1}_{W,V}(0)^\Gamma$.\hfill\break 
(i) If $\xi$ is stable or costable then the stabilizer of $\xi$ in 
$GL_\Gamma(V)$ is trivial.\hfill\break
(ii) Moreover, the stabilizer of $\xi$ in $GL_\Gamma(V)$ is trivial and the
orbit $GL_\Gamma(V)\cdot\xi$ is closed if and only if $\xi$ is stable and 
costable.
\hfill\break
(iii) The orbit $GL_\Gamma(V)\cdot\xi$ is closed if and only if there are 
subspaces
$V_1,V_2\subseteq V$ and triples $\xi_1\in\mu^{-1}_{W,V_1}(0)^{s,c,\Gamma}$, 
$\xi_2\in\mu^{-1}_{\{0\},V_2}(0)^{\Gamma}$, such that
$$V=V_1\oplus V_2,\qquad\xi=\xi_1\oplus\xi_2\and
GL_\Gamma(V_2)\cdot\xi_2\ \text{is\ closed}.$$
Moreover the splitting $\xi=\xi_1\oplus\xi_2$ is unique.
}

\vskip3mm

\noindent{\it Proof:}
Part $(i)$ is proved in \cite{N2, Lemma 3.10}.
Part $(ii)\Leftarrow$ is proved in \cite{N2, Proposition 3.24}.
Part $(ii)\Rightarrow$ is stated without proof in \cite{L, Section 2}.
Let us prove it.
Given a triple $\xi=(B,i,j)\in\mu_{W,V}^{-1}(0)^\Gamma$, let
$\xi^0=(B^0,i^0,j^0)$ be a representative of 
the unique closed orbit in $\overline{GL_\Gamma(V)\cdot\xi}$. 
Fix a one-parameter subgroup $\lambda\,:\,\CC^\times\to GL_\Gamma(V)$ such that
$$\xi^0=\lim_{t\to 0}\bigl(\lambda(t)\cdot\xi\bigr).\leqno(2.1)$$
For all $n\in\ZZ$ set
$$V^n=\{v\in V\,|\,\lambda(t)\cdot v=t^nv,\quad\forall t\}.$$
Since the limit $(2.1)$ is well-defined we have
$$\Im i\subseteq\wedge^2L\otimes\bigoplus_{n\geq 0}V^n,\qquad
\bigoplus_{n\geq 1}V^n\subseteq\Ker j,\qquad
B(V^n)\subseteq L\otimes\bigoplus_{m\geq n}V^m\quad\forall n.$$
Thus,
$$\matrix
\xi\text{\ is\ costable\ }\Rightarrow\quad V=\bigoplus_{n\geq 0}V^n\hfill\cr\cr
\xi\text{\ is\ stable\ }\Rightarrow\quad\bigoplus_{n\geq 1}V^n=\{0\}.\hfill
\endmatrix\leqno(2.2)$$
Hence, if $\xi\in\mu^{-1}_{W,V}(0)^{s,c,\Gamma}$ then $\xi=\xi^0$, i.e.
the orbit $GL_\Gamma(V)\cdot\xi$ is closed.
Part $(iii)\Rightarrow$ is proved in \cite{L, Lemma 2.30}.
The decomposition $\xi=\xi_1\oplus\xi_2$ is unique since $V_1$ 
(resp. $V_2$) is the smallest (resp. the biggest) subspace of $V$ such that
$$B(V_1)\subseteq L\otimes V_1\and\Im(i)\subseteq V_1\quad
(\text{resp.}\quad B(V_2)\subseteq L\otimes V_2\and V_2\subseteq\Ker j).$$
Let us prove Part $(iii)\Leftarrow$. Hence, fix a triple 
$\xi\in\mu_{W,V}^{-1}(0)$ admitting a splitting $\xi=\xi_2\oplus\xi_2$ as in 
$(iii)$. 
Let $\xi^0$ be a representative of the unique closed orbit in 
$\overline{GL_\Gamma(V)\cdot\xi}$.
Fix a splitting $V=V^0_1\oplus V^0_2$, $\xi^0=\xi^0_1\oplus\xi^0_2$, 
as in Part $(iii)\Rightarrow$.
Fix a one-parameter subgroup $\lambda\,:\,\CC^\times\to GL_\Gamma(V)$ such that
$\lim_{t\to 0}(\lambda(t)\cdot\xi)=\xi^0$.
Put 
$$\bar V_a^0=\lim_{t\to 0}\lambda(t)(V_a)\and
\bar\xi_a^0=\lim_{t\to 0}(\lambda(t)\cdot\xi_a)\qquad\forall a=1,2.$$

\noindent 1. First, suppose that $\lambda(t)(V_1)=V_1$ for all $t$, and
hence that $\bar V_1^0=V_1$.
The characterization of $\bar V^0_a$ and $V^0_a$ above implies that
$V_1^0\subseteq\bar V_1^0$ and $\bar V_2^0\subseteq V^0_2.$
Now, $\bar\xi^0_1\in\overline{GL_\Gamma(V_1)\cdot\xi_1}=
GL_\Gamma(V_1)\cdot\xi_1$, and thus $V_1^0=\bar V^0_1=V_1$. 
Then, a dimension count gives $V_2^0=\bar V^0_2$, and,
in particular, $V=\bar V_1^0\oplus\bar V_2^0$. 
Consider the unipotent group
$$U_\lambda=
\{u\in GL_\Gamma(V)\st\lim_{t\to 0}(\lambda(t)u\lambda(t)^{-1})=1\}.$$
Up to conjugating $\lambda$ by an element $u\in U_\lambda$ such that
$u(\bar V^0_a)=V_a$ if $a=1,2$, we can suppose that 
$\lambda=\lambda_1\times\lambda_2$
where $\lambda_a$ is a one-parameter subgroup of $GL_\Gamma(V_a)$. Then,
since $\xi_1,\xi_2$ have closed orbits, we get
$$\xi^0=\lim_{t\to 0}
\bigl(\lambda_1(t)\cdot\xi_1\oplus\lambda_2(t)\cdot\xi_2\bigr)
\,\in\, GL_\Gamma(V_1)\cdot\xi_1\oplus GL_\Gamma(V_2)\cdot\xi_2\,\subseteq\,
GL_\Gamma(V)\cdot\xi.$$

\noindent 2. Generally, fix an element $u\in U_\lambda$
such that $u(\bar V^0_1)=V_1$ and consider the one-parameter subgroup
$\lambda'(t)=u\lambda(t)u^{-1}$. We have
$$\lambda'(t)(V_1)=V_1\quad\forall t\and
\lim_{t\to 0}(\lambda'(t)\cdot\xi)=u\xi^0.$$
Thus, Part 1 implies that $u\xi^0\in GL_\Gamma(V)\cdot\xi$. Hence the orbit 
of $\xi$ is closed.
\qed

\vskip3mm

\noindent{\bf 2.3.} 
Fix $W,V\in\Rep(\Gamma)$. Fix a torus $T\subset GL_\Gamma(L)$ and put
$G_W=GL_\Gamma(W)\times T$. The obvious action of $T$ on $L$ induces an action
of $G_W$ on $M_{W,V}$ as follows 
$$(g,z)\cdot(B,i,j)=
\bigl((z\otimes 1)\circ B,(z\otimes 1)\circ i\circ g^{-1},g\circ j\bigr),
\quad\forall (g,z)\in G_W.$$
This action descends to an action of $G_W$ on the variety $\Lambda_{W,V}$.
On the other hand the action of $\Gamma$ on $L=\CC^2$ extends to $\PP^2$ 
in such a way that the line at infinity, $l_\infty$, is preserved.
Let $X_{W,V}$ be the set of isomorphism classes of
pairs $(\Ec,\phi)$ where $\Ec$ is a $\Gamma$-equivariant torsion 
free sheaf on $\PP^2$ such that
$H^1(\PP^2,\Ec(-l_\infty))\simeq V$ as a $\Gamma$-module
and $\phi$ is a $\Gamma$-invariant isomorphism
$$\phi\,:\,\Ec_{|l_\infty}{\buildrel\sim\over\to}\Oc_{l_\infty}\otimes W.$$
Put $X_W=\bigcup_V X_{W,V}$. 
Let $(e_x,e_y)$ be the canonical basis of $L$.
For any triple $\xi=(B,i,j)\in\mu_{W,V}^{-1}(0)^\Gamma$ we consider 
the following $\Gamma$-equivariant complex of sheaves on $\PP^2$ :
$$\Oc_{\PP^2}(-l_\infty)\otimes V{\buildrel a_\xi\over\longrightarrow}
\Oc_{\PP^2}\otimes\bigl((L\otimes V)\oplus W\bigr)
{\buildrel b_\xi\over\longrightarrow}
\Oc_{\PP^2}(l_\infty)\otimes\wedge^2 L\otimes V,\leqno(2.3)_\xi$$
where 
$$a_\xi=\left(\matrix zB-xe_x-ye_y\cr zj\endmatrix\right),
\quad b_\xi=\bigl(zB-xe_x-ye_y,\,zi\bigr),$$
and $x,y,z$ are homogeneous coordinates on $\PP^2$ such that 
$l_\infty=\{z=0\}.$
Let $\Hc^0_\xi$, $\Hc^1_\xi$,  and $\Hc^2_\xi$, be the cohomology sheaves 
of the complex $(2.3)_\xi$. Then (see \cite{N3, Lemma 2.6}) 

\vskip3mm

\itemitem{(a)} $\Hc^0_\xi=\{0\}$,
\itemitem{(b)} $\xi\in\mu_{W,V}^{-1}(0)^{c,\Gamma}$ if and only if the map 
$b_\xi$ is surjective,
\itemitem{(c)} $\xi\in\mu_{W,V}^{-1}(0)^{s,\Gamma}$ if and only if the map 
$a_\xi$ is fiberwise injective.

\vskip3mm

\noindent
Considering the complex $(2.3)_{\xi^*}$ simultaneously for all $\xi$ we get
a complex of sheaves on $\PP^2\times\mu_{W,V}^{-1}(0)^{s,\Gamma}$.
Since the group $GL_\Gamma(V)$ acts freely on $\mu_{W,V}^{-1}(0)^{s,\Gamma}$,
this complex descends to a complex on $\PP^2\times\Lambda_{W,V}$,
denoted by $\Cc$. The cohomology sheaves $\Hc^0(\Cc)$ and 
$\Hc^2(\Cc)$ vanish by (a), (b), and the sheaf $\Uc=\Hc^1(\Cc)$ is 
$\Gamma$-equivariant. Since $\Hc^0(\Cc)$ vanishes,
there is a $\Gamma$-invariant isomorphism
$$\Phi\,:\,\Uc_{|l_\infty\times\Lambda_{W,V}}{\buildrel\sim\over\to}
\Hc^1(\Cc_{|l_\infty\times\Lambda_{W,V}}){\buildrel\sim\over\to}
\Oc_{l_\infty\times\Lambda_{W,V}}\otimes W.$$ 
A pair $(\Uc,\Phi)$ is said to be
{\bf universal} if for any family of $\Gamma$-equivariant
torsion free sheaves $\Ec$ on $\PP^2$ parametrized by a variety $S$ 
such that

\vskip3mm

\itemitem{(d)} for all $s\in S$ we have a $\Gamma$-invariant isomorphism
$H^1\bigl(\PP^2,\Ec_{|\PP^2\times\{s\}}(-l_\infty)\bigr)\simeq V$,
\itemitem{(e)} we have a $\Gamma$-invariant trivialisation
$\phi\,:\,\Ec_{|l_\infty\times S}{\buildrel\sim\over\to}\,
\Oc_{l_\infty\times S}\otimes W,$

\vskip3mm

\noindent
there is a unique map $i_S\,:\,S\to\Lambda_{W,V}$ with 
$(i_S^*\Uc,i_S^*\Phi)\simeq(\Ec,\phi)$.
Let $X_{W,V}^0$ and $X_W^0$ be the subsets of isomorphism classes
of pairs $(\Ec,\phi)$ such that $\Ec$ is locally free. 
For any $z\in T$ let $\alpha(z)$ be the corresponding
automorphism of $\PP^2$. The group $G_W$ acts on $X_{W,V}$ in such a way that 
$$(g,z)\cdot (\Ec,\phi)=\bigl(\alpha(z)^*\Ec,g\circ\alpha(z)^*\phi\bigr),
\qquad\forall (g,z)\in G_W.$$
This action preserves the subset $X^0_{W,V}\subseteq X_{W,V}$.
For any variety $X$ let $S^nX$ be the $n$-th symmetric product and
put $SX=\bigcup_n S^nX.$

\vskip3mm

\noindent{\smc Theorem 1.} {\it There are bijections of $G_W$-sets
$$\Lambda_{W,V}{\buildrel\sim\over\to}X_{W,V}\and
N_W{\buildrel\sim\over\to} X^0_W\times (S\CC^2)^\Gamma$$
such that the map $\pi\,:\,\Lambda_W\to N_W$ is identified with the map 
$$(\Ec,\phi)\mapsto\bigl((\Ec^{**},\phi),\,\sup(\Ec^{**}/\Ec)\bigr).$$
Moreover $(\Uc,\Phi)$ is a universal pair.}

\vskip3mm

\noindent{\it Proof:} 
1. Given a triple $\xi=(B,i,j)\in\mu_{W,V}^{-1}(0)^{c,\Gamma}$, let
$\xi^0$ be a representative of 
the unique closed orbit in $\overline{GL_\Gamma(V)\cdot\xi}$. 
Fix a splitting $\xi^0=\xi^0_1\oplus\xi^0_2$, $V=V_1\oplus V_2$, as in 
Lemma 1$(iii)$, 
and fix a one-parameter subgroup $\lambda\,:\,\CC^\times\to GL_\Gamma(V)$ 
such that $\xi^0=\lim_{t\to 0}\bigl(\lambda(t)\cdot\xi\bigr).$
>From $(2.2)$ we have 
$$V=\bigoplus_{n\geq 0}V^n,\where
V^n=\{v\in V\,|\,\lambda(t)\cdot v=t^nv,\quad\forall t\}.$$
Then $\xi^0$ is the triple $\bigl(B^0,i^0,j\bigr),$
where $B^0=\lim_{t\to 0}(\lambda(t)\cdot B)$ and $i^0$ is the component
of the map $i$ in $\Lin_\Gamma(W,\wedge^2L\otimes V^0)$. Moreover,
$$\matrix
W_1=W,\hfill&\qquad V_1=V^0,\hfill&\qquad\xi^0_1=\bigl(B^0_{|V_1},i^0,j\bigr),
\hfill\cr\cr
W_2=\{0\},\hfill&\qquad V_2=\bigoplus_{n\geq 1}V^n,\hfill&
\qquad\xi^0_2=\bigl(B_{|V_2},0,0\bigr).\hfill
\endmatrix$$
Thus, we have a short exact sequence of complexes
$$0\to(2.3)_{\xi^0_2}\to(2.3)_\xi\to(2.3)_{\xi^0_1}\to 0.$$
Moreover, $\Hc^2_\xi=\{0\}$ since $\xi$ is costable and
$\Hc^1_{\xi^0_2}=\{0\}$ since $(2.3)_{\xi^0_2}$ is a Koszul complex.
Thus,
$$0\to\Hc^1_\xi\to\Hc^1_{\xi^0_1}\to\Hc^2_{\xi^0_2}\to 0\leqno(2.4)$$
is an exact sequence. Since $\xi^0_1$ is stable and costable, Claims
(b) and (c) imply that $\Hc^1_{\xi^0_1}$ is a locally free sheaf.
Since $\Hc^2_{\xi^0_2}$ has a finite length, the sheaf $\Hc^1_\xi$ is 
torsion free and its double dual is $\Hc^1_{\xi^0_1}.$
Since $\Hc^0_\xi$ vanishes we have a $\Gamma$-isomorphism
$$\phi_\xi\,:\,\Hc^1_{\xi|l_\infty}=\Hc^1\bigl((2.3)_{\xi|l_\infty}\bigr)
{\buildrel\sim\over\to}\Oc_{l_\infty}\otimes W.$$
Since $\Hc^1_\xi$ is torsion free and is trivial on $l_\infty$, we have
$H^i\bigl(\PP^2,\Hc^1_\xi(-l_\infty)\bigr)=\{0\}$ if $i=0,2$ 
(see \cite{N3, Lemma 2.3}).
Tensoring $(2.3)_\xi$ with $\Oc_{\PP^2}(-l_\infty)$ we obtain a complex
of acyclic sheaves quasi-isomorphic to $\Hc^1_\xi(-l_\infty)$.
It gives $H^1\bigl(\PP^2,\Hc^1_\xi(-l_\infty)\bigr)\simeq V.$ 
Hence we have a map
$$\mu_{W,V}^{-1}(0)^{s,\Gamma}\to X_{W,V},\quad
\xi\mapsto (\Hc^1_{\xi^*},\phi_{\xi^*}),$$
which descends to a map $\Psi\,:\,\Lambda_{W,V}\to X_{W,V}$.

\vskip3mm

\noindent 
2. It is proved in \cite{N3, Section 2.1} that $\Psi$ is surjective. The case of 
families is identical. Let us sketch it for the convenience of the reader. 
Fix an algebraic variety $S$, put $\PP=\PP^2\times S$, 
and let $\rho\,:\,\PP\to S$ be the projection.
Fix a familly $\Ec$ of $\Gamma$-equivariant torsion free sheaves on $\PP^2$
parametrized by $S$ endowed with $\Gamma$-isomorphisms 
$$\psi\,:\,R^1\rho_*\bigl(\Ec(-l_\infty)\bigr){\buildrel\sim\over\to}
\Oc_S\otimes V\and
\phi\,:\,\Ec_{|l_\infty\times S}{\buildrel\sim\over\to}\,
\Oc_{l_\infty\times S}\otimes W.$$
There is a (relative) Beilinson spectral 
sequence converging to $\Ec(-l_\infty)$ with $E_1$-term
$$E_1^{-i,j}=\Oc_{\PP^2}(-il_\infty)\boxtimes 
R^j\rho_*\bigl(\Ec\otimes\Omega^i_{\PP/S}((i-1)l_\infty)\bigr),$$
where $\Omega^i_{\PP/S}$ is the relative bundle of $i$-forms, and $i,j=0,1,2$.
Then, \cite{N3, Lemma 2.3} states that the spectral sequence degenerates at the
$E_2$ term, giving a $\Gamma$-equivariant complex
$$\Oc_{\PP^2}(-l_\infty)\boxtimes R^1\rho_*\bigl(\Ec(-2l_\infty))
{\buildrel a\over\to}
\Oc_{\PP^2}\boxtimes R^1\rho_*\bigl(\Ec\otimes\Omega^1_{\PP/S})
{\buildrel b\over\to}
\Oc_{\PP^2}(l_\infty)\boxtimes R^1\rho_*\bigl(\Ec(-l_\infty)),$$
such that $Ker\,a=Coker\,b=\{0\}$ and $Ker\,b/Im\,a=\Ec.$
Now, $\psi$ and $\phi$ give $\Gamma$-isomorphisms
$$R^1\rho_*\bigl(\Ec(-l_\infty))\simeq R^1\rho_*\bigl(\Ec(-2l_\infty))
\simeq\Oc_S\otimes V$$
$$R^1\rho_*\bigl(\Ec\otimes\Omega^1_{\PP/S})\simeq
\Oc_S\otimes(L\otimes V\oplus W),$$
and the maps $a,b,$ can be expressed as $a_\xi,b_\xi$, in $(2.3)_\xi$.
Let $(\Uc',\Phi')$ be the pull-back of $(\Uc,\Phi)$ to 
$\PP^2\times\mu_{W,V}^{-1}(0)^{s,\Gamma}$.
We have constructed a map $j_S\,:\,S\to\mu_{W,V}^{-1}(0)^{s,\Gamma}$ such that
$(\Ec,\phi)\simeq (j^*_S\Uc',j^*_S\Phi')$. If the sheaf 
$R^1\rho_*\bigl(\Ec(-l_\infty)\bigr)$ is no longer trivial but satisfies (d),
fix an open covering $(S_k)_k$ of $S$ such that
$R^1\rho_*\bigl(\Ec(-l_\infty)\bigr)_{|S_k}$ is trivial for all $k$. 
The corresponding maps
$S_k\to\mu^{-1}_{W,V}(0)^{s,\Gamma}\to\Lambda_{W,V}$ glue together.
Let $i_S$ be the resulting map.

\vskip3mm

\noindent
3. Let us prove that $\Psi$ is injective. Fix 
$\xi^1,\xi^2\in\mu^{-1}_{W,V}(0)^{s,\Gamma}$
such that $\Psi(\xi^1)\simeq\Psi(\xi^2)$. By \cite{OSS, Lemma II.4.1.3} the 
isomorphism lifts to an isomorphism of the complexes 
$(2.3)_{\xi^1}\simeq(2.3)_{\xi^2}$. This isomorphism is a triple 
$(\alpha,\beta,\gamma)\in 
GL_\Gamma(V)\times GL_\Gamma(L\otimes V\oplus W)\times GL_\Gamma(V).$
Considering the fibers at the points $[1:1:0]$, $[0:1:0]$, and $[1:0:0]$,
we get $\alpha=\gamma$ and $\beta(u\otimes v)=u\otimes\alpha(v)$ for
all $u\in L$, $v\in V$. Then, since the isomorphim is compatible with
the trivialisations $\phi_{\xi^{1*}}$ and $\phi_{\xi^{2*}}$, the map
$\beta$ fixes $W$. Hence, $\xi^2=\alpha\cdot\xi^1$. We are done.
In particular, the map $i_S$ associated to a family parametrized by $S$ 
is unique.

\vskip3mm

\noindent
4. Now, recall that points in $N_{W,V}$ are in bijection
with closed $GL_\Gamma(V)$-orbits in $\mu^{-1}_{W,V}(0)^\Gamma$.
By Claims (b), (c), we know that $\Psi$ maps 
$\mu_{W,V}^{-1}(0)^{s,c,\Gamma}/GL_\Gamma(V)$ to $X^0_{W,V}.$ 
Thus the second claim in the theorem follows from Lemma 1$(iii)$. 
By construction $\pi$ maps
the $GL_\Gamma(V)$-orbit of a triple $\xi\in\mu_{W,V}^{-1}(0)^{s,\Gamma}$ to 
the unique closed orbit in $\overline{GL_\Gamma(V)\cdot\xi}$. 
Thus the third claim follows from the exact sequence $(2.4)$. 
\qed

\vskip3mm

\noindent{\smc Remarks.}
1. Observe that (d) and (e) guarantee that the family $\Ec$ is flat over $S$, 
by \cite{H}, since all geometric fibers of $\Ec$
over $S$ have the same Hilbert polynomial. 

\noindent
2. Observe also that Part 3 of the proof implies that elements of $X_{W,V}$ 
do not have automorphisms.

\noindent
3. Finally, observe that the map $\pi$ above is the same as the one used in
\cite{B}.

\vskip3mm

\noindent{\smc Corollary.} {\it The fibers of the map 
$\pi\,:\,\Lambda_{W,V}\to N_{W,V}$ are isomorphic
to the $\Gamma$-fixpoint sets of a product of punctual $Quot$-schemes. In 
particular, $\pi$ is a proper map.}\qed

\vskip1cm

\centerline{\bf 3. The convolution algebra.} 

\vskip3mm

\noindent{\bf 3.1.}
For any linear algebraic group $G$ and any quasi-projective
$G$-variety $X$ let $K_G(X)$ and $R(G)$ be 
respectively the complexified Grothendieck group of $G$-equivariant 
coherent sheaves on $X$ and the complexified Grothendieck group
of finite dimensional $G$-modules.
If $\Ec$ is a $G$-equivariant coherent sheaf on $X$, let $\Ec$ 
denote also its class in $K_G(X)$. 
Similarly, for any $V\in\Rep(G)$, let $V$ denote also the class of $V$ in 
$R(G)$.
Recall that if $X'\subseteq X$ is a closed subvariety, then
$K_G(X')$ is identified with the complexified Grothendieck group
of the Abelian category of coherent sheaves on $X$ supported on $X'$.
Then, let $[X']\in K_G(X)$ be the class of the structural sheaf $\Oc_{X'}$ of 
$X'$.
Suppose now that $X$ is smooth and that
$\pi\,:\,X\to Y$ is a proper $G$-equivariant map, with $Y$ a 
quasi-projective (possibly singular) $G$-variety. 
Let $Z=X\times_Y X$ be the fiber product, endowed with the 
diagonal action of $G$. 
If $1\leq i,j\leq 3$, let $p_{ij}\,:\,X^3\to X^2$
be the projection along the factor not named. Obviously, 
$Z=p_{13}\bigl(p_{12}^{-1}(Z)\cap p_{23}^{-1}(Z)\bigr).$
Thus, we have the convolution map 
$$\star\,:\,K_G(Z)\times K_G(Z)\to K_G(Z),\qquad
(x,y)\mapsto Rp_{13\,*}\bigl((p_{12}^*x) \otimes^\LL(p_{23}^*y)\bigr).$$
It is known that $(K_G(Z),\star)$ is an associative algebra with unit
$[\Delta]$, where $\Delta\subseteq X\times X$ is the diagonal, see \cite{CG}.
Similarly we consider the map
$$\star\,:\,K_G(Z)\times K_G(X)\to K_G(X),\qquad
(x,y)\mapsto Rp_{1\,*}\bigl(x\otimes^\LL(p_2^*y)\bigr),$$
where $p_1$ and $p_2$ are the projections of $X\times X$ onto its factors.
Then, $K_G(X)$ is a $K_G(Z)$-module.
Now, suppose that $G$ is reductive and that $z\in G$
is a central element such that the fixpoint subvariety $X^z$ is 
compact (and smooth). Let $\iota\,:\,X^z\hookrightarrow X$ and 
$q_1\,:\,X\times X^z\to X$ be the obvious maps. 
Let $K_G(X)_z$ and $R(G)_z$ be the localized rings
with respect to the maximal ideal associated to $z$.
The direct image morphism $R\iota_*\,:\,K_G(X^z)_z\to K_G(X)_z$ is invertible
by the Localization Theorem.  Hence, since $q_1$ is proper there is a 
well-defined convolution product
$$\star\,:\,K_G(X)_z\times K_G(X)_z\to K_G(X)_z,\qquad
(x,y)\mapsto Rq_{1\,*}\bigl(x\boxtimes(R\iota_*)^{-1}(y)\bigr),$$
and the Kunneth theorem \cite{CG, Theorem 5.6.1} holds.
In particular, if $[\Delta]\in K_G(X)\otimes_{R(G)}K_G(X)$ then
$K_G(X)_z$ is a projective $R(G)_z$-module and
$$K_G(X\times X)_z\simeq K_G(X)_z\otimes_{R(G)_z}K_G(X)_z.$$
For more details on the convolution product see \cite{CG}.
We introduce two more notations.
For any equivariant vector bundle $\Ec$ and any $z\in\CC$ set 
$\Lambda_z\Ec=\sum_i(-z)^i\wedge^i\Ec$, where $\wedge^i$ stands
for the usual wedge power, and let $D\Ec$ be the determinant. 
If $z$ is a formal variable it is well-known that $\Lambda_z,D,$ induce maps
$D\,:\,K_G(X)\to K_G(X)$, $\Lambda_z\,:\,K_G(X)[[z]]\to K_G(X)[[z]]$, 
such that for all $\Ec,\Fc\in K_G(X)$ we have 
$$D(\Ec+\Fc)=D(\Ec)D(\Fc)\and
\Lambda_z(\Ec+\Fc)=\Lambda_z(\Ec)\Lambda_z(\Fc).$$
Moreover, if $\Ec,\Fc\in K_G(X)$ have the same rank then
$\Lambda_z\Ec(\Lambda_z\Fc)^{-1}$ has an expansion in $K_G(X)[[z]]$ and
$K_G(X)[[z^{-1}]]$. 
For simplicity we put $\Lambda\Ec=\Lambda_1\Ec$. If $V,W$ are
$G$-modules such that $\Lambda W$ is non zero, let
$\Lambda(V-W)$ be the class of $\Lambda V(\Lambda W)^{-1}$ in the fraction 
field of $R(G)$.

\vskip3mm

\noindent{\bf 3.2.} 
We use the convolution product in the following situation.
Consider the fiber product $Z_W=\Lambda_W\times_{N_W}\Lambda_W$.
Given $V^a\in\Rep(\Gamma)$ and 
$\xi^a=(B^a,i^a,j^a)\in\mu^{-1}_{W,V^a}(0)^\Gamma$,
$a=1,2$, put $\xi=(\xi^1,\xi^2)$ and form the complex
$$\Lin_\Gamma(V^1,V^2)\quad{\buildrel a_\xi\over\longrightarrow}
\quad
\matrix
\Lin_\Gamma(V^1,L\otimes V^2)\cr\oplus\cr 
\Lin_\Gamma(W,\wedge^2L\otimes V^2)\cr\oplus\cr\Lin_\Gamma(V^1,W)
\endmatrix
\quad
{\buildrel b_\xi\over\longrightarrow}
\quad
\Lin_\Gamma(V^1,\wedge^2L\otimes V^2),\leqno(3.1)_\xi$$
where
$$\matrix
a_\xi\,:\,f\mapsto (B^2\circ f-f\circ B^1,\,f\circ i^1,\,-j^2\circ f)
\hfill\cr\cr
b_\xi\,:\,(g,i,j)\mapsto 
B^2\wedge g+g\wedge B^1+i\circ j^1+i^2\circ j.\hfill\cr\cr
\endmatrix$$
Taking the complex $(3.1)_\xi$ simultaneously for all $\xi$
we get an equivariant complex of $G_W$-sheaves over 
$\Lambda_{W,V^1}\times\Lambda_{W,V^2}$. Suppose now that $V^1=V^2=V$.
It is proved in \cite{N2, Section 6} that the
$0$th and the $2$nd cohomology sheaves vanish.
Thus, the first cohomology sheaf, $\Hc^1$, is locally free. 
Moreover, Nakajima has constructed a section of $\Hc^1$
vanishing precisely on the diagonal 
$\Delta_{W,V}\subseteq\Lambda_{W,V}\times\Lambda_{W,V}$.
Thus
$$[\Delta_{W,V}]=\Lambda\Hc^1\in 
K_{G_W}(\Lambda_{W,V})\otimes_{R(G_W)}K_{G_W}(\Lambda_{W,V}).$$
An element $z=(q_0,t_0)\in T$ is said to be {\bf general} if for all 
$m,n\in\ZZ,$ $$q_0^mt_0^n=1\quad\Rightarrow\quad m=n=0.$$
The motivation for this technical condition is the following simple lemma.

\vskip3mm

\noindent{\smc Lemma 2.} {\it If $z$ is general, then 
$\Lambda_W^z$ is contained in $\pi^{-1}(0)$. In particular, it is compact.}

\vskip3mm

\noindent{\it Proof:} It suffices to prove that $N_{W,V}^z=\{0\}.$
Given a representative $\xi$ of a closed orbit fixed by $z$, we consider
the splitting $\xi=\xi^1+\xi^2$ introduced in Lemma $1(iii)$. Then, $\xi^2$
is identified with a point in $(S\CC^2)^\Gamma$. Since $q_0$ and $t_0$ are 
not roots of unity, necessarily $\xi^2=0$. Thus, we can suppose that
$\xi=(B,i,j)$ is stable and costable. Suppose $V\neq\{0\}$.
Fix $g\in GL_\Gamma(V)$ such that $g\cdot\xi=z\cdot\xi$.
Put $B=B_x\otimes e_x+B_y\otimes e_y$.
By costability there is a non-zero element $v\in V$ of the form
$mi(w)$ where $w\in W$ and $m$ is a word in $B_x, B_y$,
such that $B(v)=0$. Since $q_0B_x=gB_xg^{-1}$, $t_0B_y=gB_yg^{-1}$,
$q_0t_0i=g\,i$, and since $z$ is general, we have $g(v)=\lambda v$
with $\lambda\neq 1$.
Now we have also $j\,g=j$ and, thus, $j(v)=0$. 
This is in contradiction with the stability of $\xi$.
\qed

\vskip3mm

\noindent
Thus, if $z\in T$ is general by Section 3.1 we have an isomorphism
$$K_{G_W}(\Lambda_W\times\Lambda_W)_z\simeq
\End_{R(G_W)_z}K_{G_W}(\Lambda_W)_z.$$

\vskip3mm

\noindent{\bf 3.3.}
For any $k$ let $C^+_k\subseteq\Lambda_W\times\Lambda_W$
be the variety of the pairs 
$(\xi^1,\xi^2)\in\Lambda_{W,V^1}\times\Lambda_{W,V^2}$,
where $\xi^a=(B^a,i^a,j^a)$, $a=1,2$, is such that $V^2$ admits a $B^2$-stable
$\Gamma$-submodule isomorphic to $V^1$, containing $\Im (i^2)$, 
with $V^2/V^1\simeq S_k$, 
and $\xi^1$ is isomorphic to the restriction of $\xi^2$ to $V^1$. 
Let $C^-_k$ be obtained by exchanging the components of $C^+_k$. 
The varieties $C_k^\pm$ are called {\bf Hecke correspondences}.
They are closed subvarieties of $Z_W$ (see \cite{N2}). 
Let $\Vc=\mu_{W,V}^{-1}(0)^{s,\Gamma}\times_{GL_\Gamma(V)}V$
be the {\bf universal bundle} on $\Lambda_{W,V}$ and let 
$\Wc=\Oc_{\Lambda_W}\otimes W$ be the trivial sheaf. 
The sheaves $\Vc,\Wc,$ are obviously $G_W\times\Gamma$-equivariant.
Let $\Vc^a$ be the pull-back of $\Vc$ by the $a$-th projection 
$p_a\,:\,\Lambda_W\times\Lambda_W\to\Lambda_W$, and let $\Wc$ denote
also the trivial sheaf $\Oc_{\Lambda_W\times\Lambda_W}\otimes W$. 
The restriction of $\Vc^1$ to $C_k^+$ is a subbundle of $\Vc^2$
(resp. the restriction of $\Vc^2$ to $C_k^-$ is a subbundle of $\Vc^1$).
Let $\Lc^\pm_k$ be the quotient sheaf :
it is the extension by zero to $Z_W$ of an invertible sheaf on $C^\pm_k$.
Put
$$\theta=L^*-\wedge^2L^*-1\in R(G_W\times\Gamma).$$
Fix $q,t\in R(T)$ such that $L=tS^{-1}+qS$.
For any $\Gamma$-module $V$ or any $\Gamma$-equivariant sheaf $\Ec$ put
$$V_k=\Lin_\Gamma(S_k,V),\qquad v_k=\dim V_k,\and\Ec_k=\Lin_\Gamma(S_k,\Ec).$$
If $s\in\ZZ$ and $k\in\ZZ/n\ZZ$, let us define the following classes
in $K_{G_W}(Z_W)$ : 
$$\matrix
\Omega_{k,s}^+=(\Lc^+_k)^{^s}
D\bigl(\theta\Vc^2+(\wedge^2L^*)\Wc\bigr)_k^*
\hfill\cr\cr
\Omega_{k,s}^-=(\Lc^-_k)^{^{s+h_k}}\gamma_k,\hfill
\endmatrix$$
where $h_k=\dim (\theta^*V^1+W)_k$ and
$\gamma_k=q^{v^1_k-v^1_{k-1}}t^{v_k^1-v_{k+1}^1}.$
Let $\Theta_k^{^\pm}(z)$ be the expansion of
$$\Theta_k(z)=(-1)^{h_k}\gamma_k\Lambda_z\bigl((\wedge^2L^*-1)(\theta\Vc^*+\Wc^*
)\bigr)_k$$
in $K_{G_W}(Z_W)[[z^{\mp 1}]]$
(an element in $K_{G_W}(\Lambda_W)$ is identified with its direct image
in $K_{G_W}(Z_W)$ via the diagonal embedding 
$\Lambda_W\hookrightarrow Z_W$). 
More precisely, let $\Theta_{k,s}^{^\pm}\in K_{G_W}(Z_W)$ 
be such that 
$$\Theta_k^{^\pm}(z)=\sum_{s\in\NN}\Theta^{^\pm}_{k,s}\,S_k\,z^{^{\mp s}}.$$ 

\vskip3mm

\noindent{\bf 3.4.} We now consider the cyclic quiver only. Thus,
put $\gamma=\pmatrix e^{2i\pi/n} & 0\cr 0&e^{-2i\pi/n}\endpmatrix$
and $\Gamma=\left\{\gamma^k\st k\in\ZZ/n\ZZ\right\}.$
Fix a generator $S$ of $X(\Gamma)$ and set $S_k=S^{\otimes k}$ 
for all $k\in\ZZ$. We fix $T=\CC^\times\times\CC^\times$.
Let us consider the sum
$$\var_k=\sum_V(-1)^{v_k}[\Delta_{W,V}],$$ 
where $[\Delta_{W,V}]\in K_{G_W}(Z_W)$ is the class of the diagonal.
The following theorem is the main result of the paper.

\vskip3mm

\noindent{\smc Theorem 2.} {\it If $z=(q_0,t_0)$ is general 
then the map 
$$\xb^+_{k,s}\mapsto\Omega_{k,s}^{^+}\star\var_k\star\var_{k-1},\quad
\xb^-_{k,s}\mapsto\Omega_{k,s}^{^-}\star\var_{k+1}\star\var_k,\quad 
\hb^{^\pm}_{k,s}\mapsto-\Theta^{^\pm}_{k,s}\star\var_{k-1}\star\var_{k+1},$$ 
extends uniquely to an algebra homomorphism
from $\Ub_{q_0,t_0}$ to the specialization of $K_{G_W}(Z_W)$ at the point $z$.
}\qed

\vskip1cm

\centerline{\bf 4. Proof of Theorem 2.} 

\vskip3mm

\noindent In the rest of the paper $\Gamma$ and $T$ are as in Section 3.4.

\vskip3mm

\noindent{\bf 4.1.}
Fix a maximal torus $T(W)\subseteq GL_\Gamma(W)$ and set $T_W=T(W)\times T$.
Put $R_W=R(T_W)$ and $K_W=K_{T_W}(\Lambda_W)$. 
Let $\bar R_W$ be the fraction field of $R_W$ and set
$\bar K_W=\bar R_W\otimes_{R_W}K_W$.
Fix $z\in T$ and consider the fixpoint
varieties $Z_W^z,\Lambda_W^z$. 
The action by convolution of $K_{G_W}(Z_W)$, $K_{G_W}(Z_W^z)$,
on $K_{G_W}(\Lambda_W)$, $K_{G_W}(\Lambda_W^z)$, gives two 
algebra homomorphisms
$$\matrix
K_{G_W}(Z_W)\to\End_{R(G_W)}K_{G_W}(\Lambda_W)\hfill\cr\cr
K_{G_W}(Z_W^z)\to\End_{R(G_W)}K_{G_W}(\Lambda_W^z).\hfill
\endmatrix\leqno(4.1)$$
The bivariant version of the Localization Theorem \cite{CG, Theorem 5.11.10}
gives isomorphisms 
$$K_{G_W}(Z_W)_z\simeq K_{G_W}(Z_W^z)_z\and 
K_{G_W}(\Lambda_W)_z\simeq K_{G_W}(\Lambda_W^z)_z$$
which commute with the convolution product.
Thus, it identifies the maps (4.1) after localization with respect to $z$.
Suppose now that $z$ is general. 
Then $\Lambda_W^z\subseteq\pi^{-1}(0)$ and, 
thus, $Z^z_W=\Lambda_W^z\times\Lambda_W^z$.
In particular $K_{G_W}(\Lambda_W)_z$ is a faithful $K_{G_W}(Z_W)_z$-module
by Section 3.2. Observe that (see \cite{CG; Theorem 6.1.22} for instance) 
$$K_W=R_W\otimes_{R(G_W)}K_{G_W}(\Lambda_W).$$
Observe also that, since $R_W$ is free over $R(G_W)$,
$$\End_{R_W}(K_W)=R_W\otimes_{R(G_W)}\End_{R(G_W)}K_{G_W}(\Lambda_W).$$
Moreover, $\End_{R_{W,z}}K_{W,z}$ is a projective $R_{W,z}$-module 
by Sections 3.1 and 3.2.
Thus, we have a chain of injective algebra homomorphisms
$$K_{G_W}(Z_W)_z\hookrightarrow\End_{R(G_W)_z}K_{G_W}(\Lambda_W)_z
\hookrightarrow\End_{R_{W,z}}K_{W,z}\hookrightarrow\End_{\bar R_W}\bar K_W.$$
Let $x_{k,s}^\pm\in\End_{\bar R_W}\bar K_W$ be the operator such that
$$x_{k,s}^\pm(x)=Rp_{1*}\bigl(\Omega_{k,s}^\pm\otimes^\LL p_2^*(x)\bigr),
\qquad\forall x\in K_{W}.$$
The chain of embeddings above sends the image of $\Omega_{k,s}^\pm$ 
in $K_{G_W}(Z_W)_z$ to $x_{k,s}^\pm.$
Hence, it suffices to compute the relations between 
the operators $x_{k,s}^\pm$.

\vskip3mm

\noindent{\bf 4.2.} Let $\Pi$ be the set of all partitions 
$\lambda=(\lambda_1\geq\lambda_2\geq...).$ 
A partition is identified with its Ferrer diagram, whose boxes are labelled
by pairs in $\NN\times\NN$.
Let $\Ab=\CC[x,y]$ be the ring of regular functions on $L$.
The $\Gamma$-action on $L$ induces an action of $\Gamma$ on $\Ab$.
If $\lambda\in\Pi$ let $J_\lambda\subset\Ab$ be the subspace linearly spanned
by the monomials $x^iy^j$ such that $(i,j)\notin\lambda$.
The quotient space $V_\lambda=\Ab/J_\lambda$ is a $\Gamma$-module.
Let us consider the linear operators 
$B^*_\lambda\in\Lin_\Gamma(V_\lambda,L\otimes V_\lambda)$ and
$i^*_\lambda\in\Lin_\Gamma(\CC,V_\lambda)$ such that 
$$B^*_\lambda=e_x\otimes x+e_y\otimes y\and i^*_\lambda(1)=1+J_\lambda.$$ 
Put $\xi^*_\lambda=(B^*_\lambda,i^*_\lambda,0)$ and let
$\xi_\lambda=(B_\lambda,0,j_\lambda)\in\Lambda_{\CC,V_\lambda}$ 
be the dual triple (the representation $V_\lambda$ is identified with its dual).
If $W\in\Rep(\Gamma)$ is one-dimensional then 
$\Lambda_{\CC,V_\lambda}$ is isomorphic
to $\Lambda_{W,W\otimes V_\lambda}$ and, thus, $\xi_\lambda$ is identified
with a point in $\Lambda_{W,W\otimes V_\lambda}$.
The action of $T_W$ on $\Lambda_W$ has only finitely many fixpoints.
More precisely, fix $k_1,...,k_w\in\ZZ/n\ZZ$ such that
$W\simeq\bigoplus_{a=1}^wS_{k_a}$ as a $\Gamma$-module.
Then we have the following lemma.

\vskip3mm

\noindent{\smc Lemma 3.} {\it Fix $W,V\in\Rep(\Gamma)$.\hfill\break
(i) We have $\Lambda_W^{T(W)}\simeq\prod_{a=1}^w\Lambda_{S_{k_a}}$.
\hfill\break
(ii) Moreover if $W$ is one-dimensional then 
$\Lambda_{W,V}^{T}$ is the set of all the triples
$\xi_\lambda$ such that $W\otimes V_\lambda=V$ as a $\Gamma$-module.}

\vskip3mm

\noindent{\it Proof:} Part $(ii)$ is well-known, since for all 
$W\in X(\Gamma)$,
the variety $\Lambda_W$ is a subvariety of the Hilbert scheme of all
finite length subschemes in $\CC^2$ (see \cite{N3, Section 2.2}).
Fix $\xi=(B,i,j)\in\Lambda_W^{T(W)}$ and $h\in T_W$, $h$ generic. 
Fix an isomorphism $T_W\simeq(\CC^\times)^w$ and put
$h=(t_1,t_2,...,t_w)\in(\CC^\times)^w.$ 
Set $W_a=\Ker(h-t_a)\subseteq W$ for all $a=1,2,...,w$.
Let $g\in GL_\Gamma(V)$ be such that
$$(B,i\circ h^{-1},h\circ j)=\bigl((1\otimes g)\circ B\circ g^{-1},
(1\otimes g)\circ i,j\circ g^{-1}\bigr).\leqno(4.2)$$
If $a=1,2,...,w$ let $V_a\subseteq V$ be the generalized eigenspace of $g$
associated to the eigenvalue $t_a^{-1}$. Then, 
$$B(V_a)\subseteq L\otimes V_a\qquad\forall a.$$
If $\bigoplus_a V_a\neq V$ then there exists a non-zero complex number $z$
such that $\Ker(g-z^{-1})\neq\{0\}$ and $z\neq t_a$ for all $a$.
Moreover, $(4.2)$ implies that $\Ker(g-z^{-1})\subseteq\Ker j$ and is 
$B$-stable, which contradicts the stability of the triple $\xi$. Hence,
$\xi=\bigoplus_a(B_a,i_a,j_a)$ where
$$B_a\,:\,V_a\to L\otimes V_a\qquad i_a\,:\,W_a\to\wedge^2L\otimes V_a,
\qquad j_a\,:\,V_a\to W_a$$
are the restrictions of $B,i,j$.
In particular, $(B_a,i_a,j_a)$ is stable.
\qed

\vskip3mm

\noindent{\bf 4.3.}
Given a multipartition $\ul=(\lambda_1,...,\lambda_w)\in\Pi^w$, put 
$V_\ul=\bigoplus_{a=1}^wS_{k_a}\otimes V_{\lambda_a}$.
The space $V_\ul$ is endowed with the structure of a $\Gamma$-module
as in Section 4.2. Put $B_\ul=\bigoplus_{a=1}^wB_{\lambda_a}\in 
\Lin_\Gamma(V_\lambda,L\otimes V_\lambda)$ and 
$j_\ul=\bigoplus_{a=1}^w j_{\lambda_a}\in\Lin_\Gamma(W,V_\lambda)$.
Let $\xi_\ul=(\xi_{\lambda_1},...,\xi_{\lambda_w})\in\Lambda_{W,V_\ul}^{T_W}$ 
be the class of the triple $(B_\ul,0,j_\ul)$.
Let us consider the map 
$\iota_\ul\,:\,\{\bullet\}\to\Lambda_W,\,\bullet\mapsto \xi_\ul,$
and put
$$b_\ul=R\iota_{\ul\,*}(1)\in K_W,\qquad
T_\ul=\iota_\ul^*(T_{\xi_\ul}\Lambda_W)\in R_W.$$
The elements $b_\ul$, $\ul\in\Pi^w$, form a basis of the 
$\bar R_W$-vector space $\bar K_W$ by the Localization Theorem.
Let $p\,:\,\Lambda_W^{T_W}\to\{\bullet\}$ 
be the projection to a point and let
$\iota\,:\,\Lambda_W^{T_W}\to\Lambda_W$ be the embedding. 
The map $p$ is obviously proper. We consider the pairing
$$\langle\quad|\quad\rangle\,:\quad\bar K_W\otimes\bar K_W\to\bar R_W,\quad
x\otimes y\mapsto 
Rp_*\bigl((R\iota_*)^{-1}(x)\otimes^\LL(R\iota_*)^{-1}(y)\bigr).$$

\vskip3mm

\noindent{\smc Lemma 4.} {\it If $\ul,\um\in\Pi^w$ then
$\langle b_\ul|b_\um\rangle=\delta_{\ul,\um}\Lambda (T_\ul^*).$}

\vskip3mm

\noindent{\it Proof:} 
Let $q\,:\,\Lambda_W\to\{\bullet\}$ be the projection. 
The map $q$ is not proper, but the direct image of a sheaf with compact
support is well-defined. Then we have
$$\matrix
\langle b_\ul|b_\um\rangle&=
Rq_*\bigl(R\iota_{\ul *}(1)\otimes^\LL R\iota_{\um *}(1)\bigr)
\hfill\cr\cr
&=R(q\circ\iota_\ul)_*\bigl(\iota^*_\ul R\iota_{\um *}(1)\bigr),\hfill
\endmatrix$$
by the projection formula. The lemma follows since 
$\iota^*_\ul R\iota_{\um *}(1)=\delta_{\ul,\um}\Lambda (T^*_\ul).$
\qed

\vskip3mm

\noindent{\bf 4.4.}
The Hecke correspondence $C^\pm_k$ is smooth (see \cite{N2}).
If $\xi_{\um\ul}=(\xi_\um,\xi_\ul)\in C^\pm_k$ 
let $N_{\um\ul}\in R_W$ be the class of the fiber
of the normal bundle to $C^\pm_k$ in $\Lambda_W\times\Lambda_W$ at 
$\xi_{\um\ul}$ and let $\Omega^{^{(s)}}_{\um\ul}\in R_W$ be the 
restriction of $\Omega^\pm_{k,s}$ to $\xi_{\um\ul}$. 
Let us use the following notation :
$$\um\heck\ul\quad\iff\quad \xi_{\um\ul}\in C_k^+.$$ 
Then, for all $\ul\in\Pi^w$ we have
$$x_{k,s}^\pm(b_\ul)=\sum_\um
\Omega^{^{(s)}}_{\um\ul}\Lambda_{\um\ul}\,b_\um,$$
for some $\Lambda_{\um\ul}\in\bar R_W$, where
the sum is over all the multipartitions $\um$ such that
$\um\heck\ul$ (resp. $\ul\heck\um$).

\vskip3mm

\noindent{\smc Lemma 5.} {\it If $\um\heck\ul$ or $\ul\heck\um$ then 
$\Lambda_{\um\ul}=\Lambda(N^*_{\um\ul}-T_\um^*)\in\bar R_W.$}

\vskip3mm

\noindent{\it Proof:} We have 
$$\matrix
\Omega^{^{(s)}}_{\um\ul}\Lambda_{\um\ul}&=
\langle b_\um|b_\um\rangle^{-1}\langle x_{k,s}^\pm(b_\ul)|b_\um\rangle
\hfill\cr\cr
&=\langle b_\um|b_\um\rangle^{-1}
R(p\times p)_*\biggl(p_2^*R\iota_{\ul *}(1)\otimes^\LL
p_1^*R\iota_{\um *}(1)\otimes^\LL\Omega^\pm_{k,s}\biggr)\hfill\cr\cr
&=\langle b_\um|b_\um\rangle^{-1}
(\iota_\um\times\iota_\ul)^*(\Omega^\pm_{k,s})\hfill\cr\cr
&=\Omega^{^{(s)}}_{\um\ul}\Lambda(T_\um^*)^{-1}\Lambda(N^*_{\um\ul}).\hfill
\endmatrix$$
\qed

\vskip3mm

\noindent
Lemma 5 can be made more explicit as follows. 
Recall that $R(T)=\CC[q^{\pm 1}, t^{\pm 1}]$ and $L=t S^{-1}+qS.$
Put $R(T(W))=\CC[X_1^{\pm 1},...,X_w^{\pm 1}]$, in such a way that
$W=\sum_a X_aS_{k_a}$ in $R(T_W\times\Gamma)$.

\vskip3mm

\noindent{\smc Lemma 6.} 
{\it For all $\ul\in\Pi^w$ the following equality holds
$$\iota_\ul^*\Vc=\sum_a\sum_{(i,j)\in\lambda_a}q^{i}t^{j}X_aS_{k_a+i-j}.$$
}

\vskip3mm

\noindent{\it Proof:}
Recall that the fixpoint $\xi_\ul\in\Lambda_{W,V_{\ul}}^{T_W}$ is the class 
of the triple $(B_\ul,0,j_\ul)$ as in Section 4.3.
Recall also that the action of $GL_\Gamma(V_\ul)$ on 
$\mu^{-1}_{W,V_\ul}(0)^{s,\Gamma}$ is free. Let
$\rho\,:\,T_W\to GL_\Gamma(V_\ul)$ be the group homomorphism such that
for any $g=(h,z)\in T_W$,
$$\biggl((z\otimes 1)\circ B_\ul,h\circ j_\ul\biggr)=
\biggl(\bigl(1\otimes\rho(g)\bigr)\circ B_\ul\circ\rho(g)^{-1},
j_\ul\circ\rho(g)^{-1}\biggr).\leqno(4.3)$$
Thus, $V_\ul$ may be viewed as a $T_W\times\Gamma$-module.
We must prove that the class of $V_\ul$ in the Grothendieck ring is
$$\sum_a\sum_{(i,j)\in\lambda_a}q^{i}t^{j}X_aS_{k_a+i-j}.$$
If $h=(t_1,...,t_w)\in T(W)$ and $z=(q_0,t_0)\in T$ then $(4.3)$ implies that
$\rho(g)$ acts on $S_{k_a}\otimes(x^iy^j+J_\lambda)$ by the scalar
$q^i_0t^j_0t_a.$ The proof is finished.\qed

\vskip3mm

\noindent
For all $\ul\in\Pi^w$ set 
$V_\ul=\sum_a\sum_{(i,j)\in\lambda_a}q^{i}t^{j}X_aS_{k_a+i-j}$,
$$R_\ul=\sum_k\sum_{\um\heck\ul}V_{\ul/\um}\and
I_\ul=\sum_k\sum_{\ul\heck\um}V_{\um/\ul},$$
where $V_{\ul/\um}=V_\ul-V_\um$. Recall that 
$$\theta=t^{-1}S+q^{-1}S^{-1}-1-q^{-1}t^{-1}\in R(T_W\times\Gamma)$$
(see Section 3.3). For any $\ul\in\Pi^w$ put 
$H_\ul=\theta^*\, V_\ul+W\in R(T_W\times\Gamma).$

\vskip3mm

\noindent{\smc Lemma 7.} {\it For all $\ul\in\Pi^w$ we have 
$H_\ul=I_\ul-qtR_\ul.$}

\vskip3mm

\noindent{\it Proof:}
If $\xi=(B,i,j)\in\mu^{-1}_{W,V}(0)^{s,\Gamma}$, then 
the fiber of $(2.3)_\xi$ at the point $[x:y:z]=[0:0:1]$ is the 
$\Gamma$-equivariant complex
$$V{\buildrel a_\xi\over\longrightarrow}W\oplus(L\otimes V)
{\buildrel b_\xi\over\longrightarrow}
(\wedge^2L)\otimes V.$$
By Claim (2.3) (c), there is no 0-th cohomology.
Let $H^i_\ul$ be the $i$-th cohomology group corresponding to 
$\xi=\xi_\ul$, and, as usual, put $H^i_{\ul,k}=\Lin_\Gamma(S_k,H^i_\ul)$.
We have $H^1_\ul-H^2_\ul=H_\ul.$ Moreover there are isomorphisms
$$\matrix
p_2^{-1}(\xi_\ul)\cap C_k^+\simeq\PP (H_{\ul,k}^{2*}),\hfill\quad&
(\xi^1,\xi_\ul)\mapsto\wedge^2L^*\otimes(V_\ul/V^1)_k^*,\hfill\cr\cr
p_2^{-1}(\xi_\ul)\cap C_k^-\simeq\PP (H^1_{\ul,k}),\hfill\quad&
(\xi^1,\xi_\ul)\mapsto \bigl(a_{\xi^1}(V^1/V_\ul)\bigr)_k,\hfill
\endmatrix$$
(for the second map, observe that $(L\otimes V^1)_k=(L\otimes V_\ul)_k$).
Thus, $T_W\times\Gamma$ acts on $\PP(H^1_{\ul,k})$ and $\PP(H^2_{\ul,k})$
with finitely many fixpoints corresponding to the triples
$\xi_\um$ such that $\um\heck\ul$ and $\ul\heck\um$.
The lemma follows.
\qed

\vskip3mm

\noindent Recall that for any $V\in R(T_W\times\Gamma),$ 
the element $V_k\in R(T_W)$ is such that $V=\sum_kV_kS_k$.
Moreover, for any $\ul\in\Pi^w$ and any $k$ put $v_{\ul,k}=\dim V_{\ul,k}$, 
$h_{\ul,k}=\dim H_{\ul,k}$, and
$$\gamma_{\ul,k}=q^{^{v_{\ul,k}-v_{\ul,k-1}}}t^{^{v_{\ul,k}-v_{\ul,k+1}}}.$$

\vskip3mm

\noindent{\smc Lemma 8.} {\it If ${\um\heck\ul}$ then}
$$\matrix
\Omega^{^{(s)}}_{\ul\um}=(V_{\ul/\um,k})^{^{s+h_{\ul,k}}}\gamma_{\ul,k}\hfill&
\text{\it and}\quad
\Lambda_{\ul\um}=\Lambda(q^{^{-1}}t^{^{-1}}V_{\ul/\um}R^*_\um-V_{\ul/\um}I_\ul^*
)_0
\hfill\cr\cr
\Omega^{^{(s)}}_{\um\ul}=V_{\ul/\um,k}^{^s}\,D(-q^{^{-1}}t^{^{-1}}H_{\ul,k})
\hfill&\text{\it and}\quad
\Lambda_{\um\ul}=\Lambda(q^{^{-1}}t^{^{-1}}V_{\ul/\um}^*I_\ul-V_{\ul/\um}^*R_\um
)_0
.\hfill
\endmatrix$$

\vskip3mm

\noindent{\it Proof:} 
Given $V^a\in\Rep(\Gamma)$ and 
$\xi^a=(B^a,i^a,j^a)\in\mu^{-1}(0)_{W,V^a}^{\Gamma}$,
$a=1,2$, put $\xi=(\xi^1,\xi^2)$ and let us form the complex
$$\Lin_\Gamma(V^1,V^2)\quad
{\buildrel a_\xi\over\longrightarrow}
\quad
\matrix
\Lin_\Gamma(V^1,L\otimes V^2)\cr\oplus\cr 
\Lin_\Gamma(W,\wedge^2L\otimes V^2)\cr\oplus\cr\Lin_\Gamma(V^1,W)
\endmatrix
\quad
{\buildrel b_\xi\over\longrightarrow}
\quad
\matrix
\Lin_\Gamma(V^1,\wedge^2L\otimes V^2)\cr\oplus\cr\wedge^2L
\endmatrix,\leqno(4.4)_\xi$$
where
$$\matrix
a_\xi\,:\,f\mapsto (B^2\circ f-f\circ B^1,\,f\circ i^1,\,-j^2\circ f)
\hfill\cr\cr
b_\xi\,:\,(g,i,j)\mapsto 
\bigl(B^2\wedge g+g\wedge B^1+i\circ j^1+i^2\circ j,\,
tr_{V^1}(i^1\circ j)+tr_{V^2}(i\circ j^2)\bigr).\hfill\cr\cr
\endmatrix$$
Taking the complex $(4.4)_{\xi}$ simultaneously for all $\xi$
we get an equivariant complex of $T_W$-sheaves over 
$\Lambda_W\times\Lambda_W$. 
Suppose now that $V^1\subseteq V^2$ and that $V^2/V^1\simeq S_k$.
It is proved in \cite{N2, Lemma 5.2} that the
$0$th and the $2$nd cohomology sheaves vanish. 
Thus, the first cohomology sheaf, $\Hc^1$, is locally free. 
Moreover Nakajima has constructed a section of $\Hc^1$
vanishing precisely on the Hecke correspondence $C_k^+$. 
Thus, if $\um\heck\ul$ we get
$$N_{\um\ul}=\bigl(\theta^*\,V_\um^*V_\ul+qtW^*V_\ul+V_\um^*W-qt\bigr)_0.$$
Similarly we get, using the complex $(3.1)_\xi$,
$$T_\um=\bigl(\theta^*\,V_\um^*V_\um+qtW^*V_\um+V_\um^*W\bigr)_0.$$
Thus,
$$\matrix
\Lambda_{\ul\um}=\Lambda\bigl(V_{\um/\ul}H_\ul^*-q^{-1}t^{-1}\bigr)_0
\hfill\cr\cr
\Lambda_{\um\ul}=
\Lambda\bigl(q^{-1}t^{-1}V_{\ul/\um}^*H_\ul-q^{-1}t^{-1}\bigr)_0.\hfill
\endmatrix\leqno(4.5)$$
The lemma follows from Lemma 7.
\qed

\vskip3mm

\noindent{\bf 4.5.} We now prove the relations.
Let $\Theta^\pm_{\ul,k}(z)$ be the expansion of
$$\Theta_{\ul,k}(z)=(-1)^{h_{\ul,k}}\gamma_{\ul,k}
\Lambda_z\bigl((qt-1)H^*_{\ul,k}\bigr)$$
in $R(T_W)[[z^{\mp 1}]].$
Then, consider the operator $h^\pm_k(z)$ such that 
$$h_k^\pm(z)(b_\ul)=\Theta^\pm_{\ul,k}(z)\,b_\ul.$$
Put $x_k^\pm(z)=\sum_{s\in\ZZ}x_{k,s}^\pm z^{-s}$ for all $k$.

\vskip3mm

\noindent{\smc Lemma 9.} {\it The following relation holds 
$$(1-q^{^{-1}}t^{^{-1}})[x_k^+(z),x_l^-(w)]=\delta_{kl}\epsilon(z/w)
\bigl(h^+_k(z)-h^-_k(z)\bigr).$$}

\vskip3mm

\noindent{\it Proof:}  
Suppose first that $\ul,\um,\ua$ and $\ub$ are such that 
$$\um\heck\ua,\quad\ul\hecl\ua,\quad\ub\heck\ul\and\ub\hecl\um.$$
If $\ul\neq\um$ then $V_\ua+V_\ub=V_\um+V_\ul$. Thus,
$$\matrix
N^*_{\um\ua}+N_{\ul\ua}^*-T^*_\ua&=
\bigl(\theta\,V_\ub V_\ua^*+q^{-1}t^{-1}WV_\ua^*+V_\ub W^*-
2q^{-1}t^{-1}\bigr)_0\hfill\cr\cr
&=N^*_{\ub\ul}+N_{\ub\um}^*-T^*_\ub,\hfill
\endmatrix$$
and so
$\Lambda_{\um\ua}\Lambda_{\ua\ul}=\Lambda_{\um\ub}\Lambda_{\ub\ul}.$
Moreover, a direct computation gives 
$\Omega^{^{(0)}}_{\um\ua}\Omega^{^{(0)}}_{\ua\ul}=
\Omega^{^{(0)}}_{\um\ub}\Omega^{^{(0)}}_{\ub\ul}.$
Hence,
$$[x_{k,s}^+,x_{l,t}^-](b_\ul)=\delta_{kl}\,C\,b_\ul,$$
for some $C\in\bar R_W$. Let us now compute the constant $C$.
Put $k=l$ and $\um=\ul$. Lemma 8 gives
$$\matrix
\Lambda_{\ul\ub}\Lambda_{\ub\ul}&=
\Lambda\bigl(q^{-1}t^{-1}(V_{\ul/\ub}R^*_\ub+V^*_{\ul/\ub}I_\ul)-
V_{\ul/\ub}^*R_\ub-V_{\ul/\ub}I_\ul^*\bigl)_0\hfill\cr\cr
\Lambda_{\ul\ua}\Lambda_{\ua\ul}&=
\Lambda\bigl(q^{-1}t^{-1}(V_{\ua/\ul}R^*_\ul+V^*_{\ua/\ul}I_\ua)-
V_{\ua/\ul}^*R_\ul-V_{\ua/\ul}I_\ua^*\bigl)_0,\hfill
\endmatrix$$
and
$$\matrix
\Omega^{^{(0)}}_{\ul\ub}\Omega^{^{(0)}}_{\ub\ul}=
\gamma_{\ul,k}D(-q^{^{-1}}t^{^{-1}}V^*_{\ul/\ub}H_\ul)_0\hfill\cr\cr
\Omega^{^{(0)}}_{\ul\ua}\Omega^{^{(0)}}_{\ua\ul}=
\gamma_{\ua,k}D(-q^{^{-1}}t^{^{-1}}V^*_{\ua/\ul}H_\ua)_0.\hfill
\endmatrix$$
Recall that, with the notations in Section 3.1,
we have $D(\Ec)\Lambda(\Ec^*)=(-1)^{\dim \Ec}\Lambda(\Ec).$
Hence we obtain
$$[x_{k,s}^+,x_{k,t}^-](b_\ul)=(-1)^{h_{\ul,k}}\gamma_{\ul,k}
\left(-qt\sum_{\ul\heck\ua} V^{s+t}_{\ua/\ul,k}
{\Lambda\bigl(q^{^{-1}}t^{^{-1}}V_{\ua/\ul}R^*_\ul+qtV_{\ua/\ul}I^*_\ua\bigr)_0
\over
\Lambda\bigl(V_{\ua/\ul}R^*_\ul+V_{\ua/\ul}I^*_\ua\bigr)_0}\right.+$$
$$+\left.\sum_{\ub\heck\ul}V^{s+t}_{\ul/\ub,k}
{\Lambda\bigl(q^{^{-1}}t^{^{-1}}V_{\ul/\ub}R^*_\ub+qtV_{\ul/\ub}I^*_\ul\bigr)_0
\over
\Lambda\bigl(V_{\ul/\ub}R^*_\ub+V_{\ul/\ub}I^*_\ul\bigr)_0}\right)b_\ul.
$$
Hence, Lemma 9 follows from Lemma 7 and the following fact.

\vskip3mm

\noindent{\smc Fact.} {\it Let $q$ and $a_i$, $i\in I$, be formal variables. 
Fix a partition $I=I_1\coprod I_2$. For all $s\in\ZZ$ 
consider the following 1-form  
$$\omega_s=z^{s-1}\biggl(\prod_{j\in I_1}{1-q^{^{-1}}z/a_j\over 1-z/a_j}\biggr)
\biggl(\prod_{j\in I_2}{1-qz/a_j\over 1-z/a_j}\biggr)dz.$$
Then,
$${\res_0\omega_s+\res_\infty\omega_s\over 1-q^{^{-1}}}=-q\sum_{i\in I_2}a_i^s
\biggl(\prod_{j\in I_1}{1-q^{^{-1}}a_i/a_j\over 1-a_i/a_j}\biggr)
\biggl(\prod_{j\in I_2\setminus\{i\}}{1-qa_i/a_j\over 1-a_i/a_j}\biggr)+$$
$$+\sum_{i\in I_1}a_i^s
\biggl(\prod_{j\in I_1\setminus\{i\}}{1-q^{^{-1}}a_i/a_j\over 1-a_i/a_j}\biggr)
\biggl(\prod_{j\in I_2}{1-qa_i/a_j\over 1-a_i/a_j}\biggr)
.$$
}
\qed

\vskip3mm

\noindent{\smc Lemma 10.} {\it We have
$$\matrix
(w-qtz)^{\pm 1}x_k^\pm(w)x_k^\pm(z)=(qtw-z)^{\pm 1}x_k^\pm(z)x_k^\pm(w)
\hfill\cr\cr
(tw-z)\,x_{k+1}^+(w)x_k^+(z)=(qz-w)\,x_k^+(z)x_{k+1}^+(w)
\hfill\cr\cr
(w-tz)\,x_{k+1}^-(w)x_k^-(z)=(z-qw)\,x_k^-(z)x_{k+1}^-(w)
\hfill\cr\cr
x_l^\pm(w)x_k^\pm(z)=x_k^\pm(z)x_l^\pm(w)\hfill\text{if}\quad l\neq k,k\pm 1.
\endmatrix$$
}

\vskip3mm

\noindent{\it Proof:} 
For all $k,l$ set
$$A_\un=\{\um\in\Pi^w\st\un\hecl\um\heck\ul\}\and
A'_\un=\{\um'\in\Pi^w\st\un\heck\um'\hecl\ul\}.$$
Then we have
$$\matrix
x_l^+(w)x_k^+(z)(b_\ul)=\sum_\un\sum_{\um\in A_\un}
\epsilon(w^{-1}V_{\um/\un})\,\epsilon(z^{-1}V_{\ul/\um})\,
\Omega^\um\Lambda^\um\,b_\un\hfill\cr\cr
x_k^+(z)x_l^+(w)(b_\ul)=\sum_\un\sum_{\um'\in A'_\un}
\epsilon(w^{-1}V_{\ul/\um'})\,\epsilon(z^{-1}V_{\um'/\un})\,
\Omega^{\um'}\Lambda^{\um'}\,b_\un,\hfill
\endmatrix$$
where $\Lambda^{\um}=\Lambda_{\un\um}\Lambda_{\um\ul}$
and $\Omega^{\um}=\Omega^{^{(0)}}_{\un\um}\Omega^{^{(0)}}_{\um\ul}.$
In particular if we have
$$V_{\ul/\um}=V_{\um'/\un}\and V_{\um/\un}=V_{\ul/\um'},\leqno(4.6)$$
then
$$\matrix
\Lambda^{\um}=\Lambda\bigl(\theta\,V_{\um/\un}V^*_{\ul/\um}-
\theta\,V^*_{\um/\un}V_{\ul/\um}\bigr)_0\,\Lambda^{\um'}\hfill\cr\cr
\Omega^{\um}=D(\theta\,V_{\un/\um})_k D(\theta\,V_{\ul/\um})_l\,
\Omega^{\um'}.\hfill
\endmatrix\leqno(4.7)$$
Suppose first that $l\neq k\pm 1$. Then, there is a bijection 
$A_\un{\buildrel\sim\over\to}A'_\un,$ $\um\mapsto\um',$
such that (4.6) holds. If $l=k$ then (4.7) gives
$$(1-V_{\ul/\um}^*V_{\um/\un})(1-q^{-1}t^{-1}V_{\ul/\um}^*V_{\um/\un})
\Lambda^{\um}
=(1-V_{\ul/\um}V_{\um/\un}^*)(1-q^{-1}t^{-1}V_{\ul/\um}V_{\um/\un}^*)
\Lambda^{\um'},$$
and 
$V_{\um/\un}^{-2}\,\Omega^{\um}=V_{\ul/\um}^{-2}\,\Omega^{\um'}.$
Thus 
$$(w-qtz)\,\epsilon(w^{-1}V_{\um/\un})\,
\epsilon(z^{-1}V_{\ul/\um})\,\Omega^{\um}\Lambda^{\um}=$$
$$=(qtw-z)\,\epsilon(w^{-1}V_{\ul/\um'})\,
\epsilon(z^{-1}V_{\um'/\un})\,\Omega^{\um'}\Lambda^{\um'}.$$
If $l\neq k,k\pm 1$ then $(4.7)$ gives
$\Lambda^{\um}=\Lambda^{\um'}$, $\Omega^{\um}=\Omega^{\um'}$, and 
the first claim of the lemma follows.
Suppose now that $l=k+1$. Set
$$\matrix
B_\un=A_\un\setminus
\{\um\in A_\un\st V_{\ul/\um}=tS^{-1}V_{\um/\un}\}
\hfill\cr\cr
B'_\un=A'_\un\setminus
\{\um'\in A'_\un\st V_{\ul/\um'}=qSV_{\um'/\un}\}.
\hfill
\endmatrix$$
Then,
$$\matrix
(t-z/w)\biggl(x_{k+1}^+(w)x_k^+(z)(b_\ul)-\sum_\un\sum_{\um\in B_\un}
\epsilon(w^{-1}V_{\um/\un})\,\epsilon(z^{-1}V_{\ul/\um})\,
\Omega^{\um}\Lambda^{\um}\,b_\un\biggr)=0\hfill\cr\cr
(q-w/z)\biggl(x_k^+(z)x_{k+1}^+(w)(b_\ul)-\sum_\un\sum_{\um'\in B'_\un}
\epsilon(w^{-1}V_{\ul/\um'})\,\epsilon(z^{-1}V_{\um'/\un})\,
\Omega^{\um'}\Lambda^{\um'}\,b_\un\biggr)=0.\hfill
\endmatrix$$
Now, there is a bijection 
$B_\un{\buildrel\sim\over\to}B'_\un,$ $\um\mapsto\um',$
such that (4.6) holds. Then, (4.7) implies that if $\um\in B_\un$ then
$$(1-t^{-1}V_{\um/\un,l}^*V_{\ul/\um,k})\,\Lambda^{\um}=
(1-q^{-1}V_{\um/\un,l}V_{\ul/\um,k}^*)\,\Lambda^{\um'},$$
and
$q^{-1}V_{\um/\un,l}\,\Omega^{\um}=t^{-1}V_{\ul/\um,k}\,\Omega^{\um'}.$
The formulas for the operators $x_k^-$ are proved in a similar way.
It suffices to observe that if $\um\in A_\un$, $\um'\in A'_\un$, and (4.6) 
holds, then we have the following identities which are very similar to (4.7) :
$$\matrix
\Lambda_{\ul\um}\Lambda_{\um\un}=\Lambda\bigl(\theta\,V_{\um/\un}V^*_{\ul/\um}-
\theta\,V^*_{\um/\un}V_{\ul/\um}\bigr)_0\,\Lambda_{\ul\um'}\Lambda_{\um'\un}
\hfill
\cr\cr
\Omega^{^{(0)}}_{\ul\um}\Omega^{^{(0)}}_{\um\un}=D(\theta V_{\un/\um})_k 
D(\theta V_{\ul/\um})_l
\Omega^{^{(0)}}_{\ul\um'}\Omega^{^{(0)}}_{\um'\un}.\hfill
\endmatrix$$
\qed

\vskip3mm

\noindent{\smc Lemma 11.} {\it We have
$$\matrix
h_l^+(w)x_k^+(z)=x_k^+(z)h_l^+(w)\hfill\text{if}\quad l\neq k,k\pm 1\cr\cr
(w-qtz)\,h_k^+(w)x_k^+(z)=(qtw-z)\,x_k^+(z)h_k^+(w)\hfill\cr\cr
(tw-z)\,h_{k+1}^+(w)x_k^+(z)=(qz-w)\,x_k^+(z)h_{k+1}^+(w).\hfill
\endmatrix$$
}

\vskip3mm

\noindent{\it Proof:} For all $\ul,\um\in\Pi^w$ we have
$$\Theta_{\ul,l}(w)=(-1)^{h_{\ul,l}+h_{\um,l}}\gamma_{\ul,l}\gamma^{^{-1}}_{\um,
l}
\Lambda_w\bigl((\theta^*-\theta)V^*_{\ul/\um}\bigr)_{-l}\Theta_{\um,l}(w).$$
Thus,
$$h_k^+(w)\,x_k^+(z)\,(1-wz^{-1})(qt-wz^{-1})=
(1-wz^{-1})(1-qtwz^{-1})\,x_k^+(z)\,h_k^+(w),$$
and
$$(tw-z)\,h_{k+1}^+(w)\,x_k^+(z)=(qz-w)\, x_k^+(z)\,h_{k+1}^+(w).$$
\qed

\vskip3mm

\noindent{\smc Lemma 12.} {\it We have
$$t^{^{\pm 1}}\, x_k^+(z_{_1}) x_k^+(z_{_2}) x_{k\pm 1}^+(w)+
\bigl((qt)^{^{\pm 1}}+1\bigr) x^+_k(z_{_1}) x_{k\pm 1}^+(w) x_k^+(z_{_{2}})+$$
$$+q^{^{\pm 1}}\, x_{k\pm 1}^+(w) x_k^+(z_{_{1}}) x_k^+(z_{_2})+
\{z_{_1}\leftrightarrow z_{_2}\}=0.$$
}

\vskip3mm

\noindent{\it Proof:}
Fix $l=k+1$, fix $\ul,\uo\in\Pi^w$, and consider the following sets
$$\matrix
A=\{\um,\un\in\Pi^w\st\uo\hecl\un\heck\um\heck\ul\}\hfill\cr\cr
A'=\{\um',\un'\in\Pi^w\st\uo\heck\un'\hecl\um'\heck\ul\}\hfill\cr\cr
A''=\{\um'',\un''\in\Pi^w\st\uo\heck\un''\heck\um''\hecl\ul\}.\hfill
\endmatrix$$

\vskip3mm

\noindent
1. First, suppose that $A, A'$ and $A''$ are nonempty.
Then, there are bijections
$$\matrix
A{\buildrel\sim\over\to}A',\quad (\um,\un)\mapsto(\um',\un')\hfill\cr\cr
A{\buildrel\sim\over\to}A'',\quad (\um,\un)\mapsto(\um'',\un'').\hfill
\endmatrix$$
such that
$$V_{\un/\uo}=V_{\um'/\un'},\quad V_{\um/\un}=V_{\un'/\uo},\quad
V_{\ul/\um}=V_{\ul/\um'},$$
$$V_{\un/\uo}=V_{\ul/\um''},\quad V_{\um/\un}=V_{\un''/\uo},\quad
V_{\ul/\um}=V_{\um''/\un''}.$$
Put
$$\Lambda^{\un\um}=\Lambda_{\uo\un}\Lambda_{\un\um}\Lambda_{\um\ul},\and
\Omega^{\un\um}=\Omega^{^{(0)}}_{\uo\un}\Omega^{^{(0)}}_{\un\um}\Omega^{^{(0)}}_
{\um\ul}.$$
Then, (4.5) gives
$$\matrix
\Lambda^{\un\um}&=\Lambda^{\un''\um''}\Lambda
\bigr(\theta V_{\ul/\un}^*V_{\un/\uo}-\theta V_{\ul/\un}V_{\un/\uo}^*\bigl)_0
\hfill\cr\cr
&=\Lambda^{\un'\um'}\Lambda
\bigr(\theta V_{\um/\un}^*V_{\un/\uo}-\theta V_{\um/\un}V_{\un/\uo}^*\bigl)_0
\hfill\cr\cr
\Omega^{\un\um}&=\Omega^{\un''\um''}
\Bigl(D(-\theta V_{\un/\uo})_k\Bigr)^2D(-\theta V_{\un/\ul})_l
\hfill\cr\cr
&=\Omega^{\un'\um'}D(-\theta V_{\un/\uo})_kD(-\theta V_{\un/\um})_l.
\hfill
\endmatrix$$
Thus,
$$\matrix
\Omega^{\un\um}\Lambda^{\un\um}&=q^2\Lambda(
q^{-1}V^*_{\ul/\un,k}V_{\un/\uo,l}-tV^*_{\ul/\un,k}V_{\un/\uo,l})
\Omega^{\un''\um''}\Lambda^{\un''\um''}\hfill\cr\cr
&=-q\Lambda(
q^{-1}V^*_{\um/\un,k}V_{\un/\uo,l}-tV^*_{\um/\un,k}V_{\un/\uo,l})
\Omega^{\un'\um'}\Lambda^{\un'\um'}.\hfill
\endmatrix$$
Put
$$\matrix
E(z_1,z_2)&=\epsilon(w^{-1}V_{\un/\uo})\epsilon(z_1^{-1}V_{\um/\un})
\epsilon(z_2^{-1}V_{\ul/\um})\Omega^{\un\um}\Lambda^{\un\um}\hfill\cr\cr
E'(z_1,z_2)&=\epsilon(w^{-1}V_{\um'/\un'})\epsilon(z_1^{-1}V_{\un'/\uo})
\epsilon(z_2^{-1}V_{\ul/\um'})\Omega^{\un'\um'}\Lambda^{\un'\um'}\hfill\cr\cr
E''(z_1,z_2)&=\epsilon(w^{-1}V_{\ul/\um''})\epsilon(z_1^{-1}V_{\um''/\un''})
\epsilon(z_2^{-1}V_{\un''/\uo})\Omega^{\un''\um''}\Lambda^{\un''\um''}.\hfill
\endmatrix$$
Then,
$$\matrix
E'(z_1,z_2)&={\ds -q^{-1}{z_1-tw\over z_1-q^{-1}w}E(z_1,z_2)}\hfill\cr\cr
E''(z_1,z_2)&={\ds
q^{-2}{(z_1-tw)(z_2-tw)\over (z_1-q^{-1}w)(z_2-q^{-1}w)}E(z_1,z_2).}\hfill
\endmatrix$$
Using (4.5) it is easy to see that the expression $E(z_1,z_2)$ has the form
$$E(z_1,z_2)={z_1-q^{-1}t^{-1}z_2\over z_1-z_2}\,S,$$
where the factor $S$ is symmetric in $z_1,z_2$. Moreover, observe that
$$qt+{(z_1-tw)(z_2-tw)\over(z_1-q^{-1}w)(z_2-q^{-1}w)}=$$
$$={(z_1-tw)(z_2-qtz_1)\over(z_1-q^{-1}w)(z_2-z_1)}+
{(z_2-tw)(z_1-qtz_2)\over(z_2-q^{-1}w)(z_1-z_2)}.$$
Thus we finally get
$$E(z_1,z_2)+(q+t^{-1})E'(z_1,z_2)+qt^{-1}E''(z_1,z_2)+\{z_1\leftrightarrow 
z_2\}=0.$$

\vskip3mm

\noindent
2. If $A, A'$ or $A''$ is the empty set then either
$A'\neq\emptyset$ or $A=A'=A''=\emptyset.$
Moreover, if $A'\neq\emptyset$ then, either $A\simeq A'$ and $A''=\emptyset$,
or $A\simeq A''$ and $A'=\emptyset$. Then proceed as in part 1.

\qed

\vskip1cm

\noindent{\it Acknowledgements.}
{\eightpoint{Part of this work was done while the second author was visiting
the Institute for Advanced Study at Princeton. The second author
is grateful to G. Lusztig for his kind invitation. We are also
grateful to V. Ginzburg for his interest and his encouragements.}}

\vskip3cm
{\eightpoint{
$$\matrix\format\l&\l&\l&\l\\
\phantom{.} & {\text{Michela Varagnolo}}\phantom{xxxxxxxxxxxxx} &
{\text{Eric Vasserot}}\\
\phantom{.}&{\text{D\'epartement de Math\'ematiques}}\phantom{xxxxxxxxxxxxx} &
{\text{D\'epartement de Math\'ematiques}}\\
\phantom{.}&{\text{Universit\'e de Cergy-Pontoise}}\phantom{xxxxxxxxxxxxx} &
{\text{Universit\'e de Cergy-Pontoise}}\\
\phantom{.}&{\text{2 Av. A. Chauvin}}\phantom{xxxxxxxxxxxxx} & 
{\text{2 Av. A. Chauvin}}\\
\phantom{.}&{\text{95302 Cergy-Pontoise Cedex}}\phantom{xxxxxxxxxxxxx} & 
{\text{95302 Cergy-Pontoise Cedex}}\\
\phantom{.}&{\text{France}}\phantom{xxxxxxxxxxxxx} & 
{\roman{France}}\\
&{\text{email: varagnol\@math.pst.u-cergy.fr}}\phantom{xxxxxxxxxxxxx} &
{\text{email: vasserot\@math.pst.u-cergy.fr}}
\endmatrix$$
}}
\enddocument